\documentclass[12pt,a4paper]{amsart}
\usepackage{geometry}                
\usepackage{amssymb}
\usepackage{epstopdf}
\usepackage{mathtools}
\usepackage[ruled,boxed,commentsnumbered,norelsize]{algorithm2e}
\usepackage{hyperref}
\usepackage{anysize}
\marginsize{3cm}{3cm}{3cm}{3cm}
\renewcommand{\baselinestretch}{1.25}
\newtheorem{lemma}{Lemma}
\newtheorem{proposition}{Proposition}
\newtheorem{theorem}{Theorem}
\newtheorem{remark}{Remark}
\newtheorem{corollary}{Corollary}

\newtheorem{problem}{Problem}

\theoremstyle{definition}
\newtheorem{definition}{Definition}

\newcounter{enu}
\newcounter{enubis}

\newenvironment{mylist}
{\begin{list}{P--\arabic{enu}}{\usecounter{enu}}\setcounter{enu}{\theenubis}} 
{\setcounter{enubis}{\theenu}\end{list}}

\DeclareMathOperator{\nm}{N} 
\newcommand{\Q}{\mathord{\mathbb Q}}
\newcommand{\G}{\mathord{\mathbb G}}
\newcommand{\F}{\mathord{\mathbb F}}
\newcommand{\Z}{\mathord{\mathbb Z}}
\newcommand{\N}{\mathord{\mathbb N}}

\title[{Quadratic form representations via generalized continuants}
]{Quadratic form representations via generalized continuants}

\author{Charles Delorme}

\email{\texttt{cd@lri.fr}}

\author{Guillermo Pineda-Villavicencio}
\address{Centre for Informatics and Applied Optimisation, Federation University Australia}

\email{\texttt{work@guillermo.com.au}}



\begin{document}

\begin{abstract}
H.~J.~S. Smith proved Fermat's two-square theorem using the notion of palindromic continuants. In this paper we extend Smith's approach
to proper binary quadratic form representations in some commutative Euclidean rings, including rings of integers and rings
of polynomials over fields of odd characteristic. Also, we present new deterministic algorithms for finding the corresponding proper representations.
\end{abstract}

\maketitle
\section{Introduction}
 Fermat's two-square theorem is without doubt a remarkable result. Many proofs of the theorem have been provided; see, for instance, \cite{Zag90,Her1848,Ser1848,CELV99,Bri72}. It is also true that most proofs have much in common, for instance, Smith's proof is very similar to Hermite's \cite{Her1848}, Serret's \cite{Ser1848}, and Brillhart's \cite{Bri72}.

Let us recall here, for convenience, suitable definitions of Euclidean rings and continuant.
\begin{definition}[{\cite[Section~2.15]{Jacobson}}]
\emph{Euclidean rings} are rings $R$ with no zero divisors which are endowed with a
Euclidean function $\nm$  from $R$ to  the nonnegative integers such that
for all $\tau_1, \tau_2\in R$ with $\tau_1\ne 0$, there exist $q, r \in R$
such that $\tau_2 = q\tau_1 + r$ and $\nm(r) < \nm(\tau_1)$.
\end{definition}

Well-known examples of Euclidean rings include the integers with $\nm(u)=|u|$, and the polynomials over
a field  with $\nm(0)=0$ and $\nm(P)=2^{\mathrm{degree}(P)}$. In this paper we only  consider Euclidean commutative rings.

\begin{definition}[Continuants in arbitrary rings, {\cite[Sec.~6.7]{GKP}}]\label{def:continuant}
Let $Q$ be  a sequence  of elements $q_1,q_2,\ldots,q_n$ of a ring $R$.
We associate with  $Q$ an element $[Q]$ of $R$ via the following recurrence formula
\begin{align*}
[\,]&=1, [q_1]=q_1, [q_1,q_2]=q_1q_2+1, \text{~and}\\
[q_1,q_2,\ldots,q_n]&=[q_1,\ldots,q_{n-1}]q_n+[q_1,\ldots,q_{n-2}], \text{~ if~} n\ge 3.
\end{align*}
The value $[Q]$ is called the \emph{continuant} of the sequence
$Q$.
\end{definition}

Properties of continuants in commutative rings are given by Graham et al.~\cite[Sec.~6.7]{GKP}.

\begin{lemma}[Carroll's identity, {\cite{Dod1866}}] \label{lemm:LC_Ident}
Let $C$ be an $n\times n$ matrix in a commutative ring. Let $C_{i_1,\ldots, i_s;j_1,\ldots, j_s}$
denote the matrix obtained from $C$ by omitting the rows $i_1,\ldots, i_s$ and the columns
$j_1,\ldots, j_s$. Then
\[
\det(C)\det(C_{i,j;i,j})=\det(C_{i;i})\det(C_{j;j})-\det(C_{i;j})\det(C_{j;i})
\]
where $\det(M)$ denotes the determinant of a matrix $M$, and the determinant of the
$0\times0$ matrix is 1 for convenience.
\end{lemma}

The use of Carroll's identity provides two more properties.
 
\begin{mylist}
\item\label{lewiscarroll}
 $[q_1,q_2,\ldots, q_{n}][q_2,\ldots , q_{n-1}]= [q_1,\ldots,q_{n-1}][q_2,\ldots , q_n]+(-1)^{n}$
($n\ge 2$).
\item\label{PPreverse}
$[q_1,q_2,\ldots,q_n]=[q_n,\ldots,q_2,q_1]$.
\end{mylist}

Given two elements $m_1$ and $m_2$ in a Euclidean ring $R$, the Euclidean algorithm outputs a sequence $(q_1,q_2,\ldots,q_n)$ of quotients and a greatest common divisor (gcd) $h$ of $m_1$ and $m_2$. A sequence of quotients given by the Euclidean algorithm is called a
\emph{continuant representation} of $m_1$ and $m_2$ as we have the equalities
$m_1=[q_1,q_2,\ldots,q_n]h$ and
$m_2=[q_2,\ldots,q_n] h$, unless $m_2=0$.

A representation of an element $m$ by the form $Q(x,y)=\alpha x^2+ \gamma xy+ \beta y^2$ is called {\it proper} if
$\gcd(x,y)=1$. In this paper we are mostly concerned with proper representations.

For us the quadratic forms $f(x,y)=ax^2+bxy+c y^2$ and $g(x,y)=Ax^2+Bxy+Cy^2$ are {\it equivalent} if there is a $2\times 2$ matrix $M=(a_{ij})$ with determinant 1 such that $g(x,y)=f(a_{11}x+a_{12}y, a_{21}x+a_{22}y)$. For equivalent forms $f$ and $g$ it follows that an element $m$ is (properly) represented by $f$ iff $m$ is (properly) represented by $g$.

\subsection{Our work}

This paper can be considered as a follow-up to our earlier paper \cite{DP11}. In that paper we studied the use of continuants in some integer representations
(e.g., sums of four squares) and sums of two squares in rings of polynomials over fields of characteristic different from 2.  Here we deal with the following problems. We let $u$ denote a unit in a ring; the ring under consideration will become clear from the context.

 \begin{problem}[From $x^2+gxy+hy^2$ to $z^2+gz+h$] If $m=u(x^2+gxy+hy^2)$ and $x,y$ are coprime, can we find $z$ such that $z^2+gz+h$ is a multiple of $m$ using ``continuants''?
\end{problem}
\begin{problem}[From $z^2+gz+h$ to $x^2+gxy+hy^2$] \label{prob:2}If $m$ divides $z^2+gz+h$, can we find $x,y$ such that $m=u(x^2+gxy+hy^2)$ using ``continuants''?
\end{problem}

We emphasize that, while Problem \ref{prob:2} has a positive answer in some situations (see below), in general it has a negative answer. More information  is given in Subsection \ref{subsec:Q(z,1)ToQ(x,y)}.
 
In this paper we use a generalization of  continuants \cite{She72} to produce proper representations $Q(x,y)=x^2+ g xy+ h y^2$, up to multiplication by a unit $u$,  of an element $m$ in some Euclidean rings. This generalization allows us to present the following new deterministic algorithms. 
\begin{enumerate}
\item Algorithm \ref{DivAlg_FromQ(x,y)ToQ(z,1)}: for every $m$ in a commutative Euclidean ring, it finds a solution $z_0$ of $Q(z,1)\equiv0\pmod{m}$, given a proper representation $uQ(x,y)$ of $m$.
\item Algorithm \ref{PolyDivAlg}: for every polynomial $m\in \mathbb{F}[X]$, where $\mathbb{F}$ is a field of odd characteristic, it finds a proper representation $u(x^2+hy^2)$ of $m$, given a solution $z_0$ of $Q(z,1)\equiv0\pmod{m}$. Here $h$ is a polynomial in $\mathbb{F}[X]$ of degree at most one.
\item Algorithm \ref{DivAlgI}: for all negative fundamental discriminants of class number one, it finds a representation $uQ(x,y)$ of an integer $m$, given a solution $z_0$ of $Q(z,1)\equiv0\pmod{m}$. 
\end{enumerate}

A simple modification of Algorithm \ref{DivAlgI} produces representations $uQ(x,y)$ for some positive discriminants of class number one, including all the determinants studied by Matthews \cite{Mat02}. This modification is discussed in Section \ref{sec:final}.
 
 Recall the {\it class number} of a determinant $\Delta\in \mathbb Z$ \cite[p.~7]{bue89} gives  the number of equivalence classes of integral binary quadratic forms with discriminant $\Delta$. As is customary,  we ignore negative definite forms; see \cite[p.~7]{bue89} and \cite[p.~152]{NivZucMon91}.   
 
When trying to extend Smith's approach to other Euclidean rings $R$, one is confronted by the lack of uniqueness of the continuant representation. The uniqueness of the
continuant representation boils down to the uniqueness of the quotients
and the remainders in the division algorithm. This uniqueness is achieved
only when $R$ is a field or $R = \F[X]$, the polynomial algebra over a field $\F$
\cite{Jod67} (considering the degree as the Euclidean function).

\subsection{A short review of related results}
Let $p$ be a prime  number of the form $4k+1$. In his proof of Fermat's two-square theorem,
Smith \cite{CELV99} first shows the existence of a palindromic sequence
$Q=(q_1,\ldots , q_s,q_s,\ldots ,q_1)$ such that $p=[Q]$ through an elegant  parity argument. This sequence then allows him to derive a solution for $z^2+1\equiv0 \pmod{p}$ and a representation $x^2+y^2$ for $p$.

With regard to the question of finding square roots modulo a prime $p$, Schoof \cite{Sch1985} presented a deterministic algorithm and Wagon \cite{Wag90} wrote an interesting article on the topic.

Brillhart's refinement \cite{Bri72} of Smith's construction took full advantage
of the palindromic structure of the sequence $(q_1,\dots,q_{s-1},q_s,q_s,q_{s-1},\dots,q_1)$
given by the Euclidean algorithm on $p$ and $z_0$, a solution of $z^2+1\equiv0\pmod{p}$. He noted that the Euclidean
algorithm gives the remainders
\begin{align*}
r_i&=[q_{i+2},\dots,q_{s-1},q_s,q_s,q_{s-1},\dots,q_1] \text{ ($i=1,\ldots,2s-1$), and}\\
r_{2s}&=0.
\end{align*}
So, by virtue of Smith's construction, rather than computing the whole sequence, we only
need to obtain $x=r_{s-1}=[q_s,q_{s-1},\ldots,q_1]$ and $y=r_s=[q_{s-1},\ldots,q_1]$. In this case, we have $y<x<\sqrt{p}$, Brillhart's stopping criterium.

In the ring of integers, Cornacchia \cite{Cor1908} extended Smith's ideas to cover representations of numbers $m=p$, with $p$ prime, by forms $x^2+hy^2$. It has been noticed that Cornacchia's algorithm can be used to obtain representations for all $1\le h<m$, with $m$ not necessarily prime \cite{Mag04}. Further extensions of Smith's and Brillhart's ideas have appeared in the literature \cite{HMW90,Wil1994,Wil1995}, where the authors
provided algorithms for finding proper representations of natural numbers as primitive, positive-definite, integral and binary quadratic forms. Matthews \cite{Mat02} provided representations of certain integers as $x^2-hy^2$, where $h = 2, 3, 5, 6, 7$. In all these papers continuants have featured as numerators (and denominators) of continued fractions. For instance, the continuant $[q_1,q_2,q_3]$ equals the numerator of the continued fraction $q_1 + \cfrac{1}{q_2+\cfrac{1}{q_3}}$, while the continuant $[q_2,q_3]$ equals its denominator.

For us the set of natural numbers $\N$ includes the zero.

Concerning other rings, one of the most important results is due to Choi, Lam, Reznick and Rosenberg \cite{CLR1980}.  They \cite{CLR1980} proved the following theorem.
\begin{theorem}[{\cite[Thm.~2.5]{CLR1980}}]
 Let $R$ be an integral domain, let $\mathbb{F}_R$ be its field of fractions, let $-h$ be a non-square in $\mathbb{F}_R$, and let $R[\sqrt{-h}]$ be the smallest ring containing $R$ and $\sqrt{-h}$. 
 
 If both $R$ and $R[\sqrt{-h}]$ are UFDs (unique factorisation domains), then the following assertions hold.
 
\begin{itemize}
\item[(1)] Any element $m\in R$ which is representable by the form $x'^2+hy'^2$ with $x',y'\in \mathbb{F}_R$ is also representable by the form $x^2+hy^2$ with $x,y\in R$.

\item[(2)] Any element $m\in R$ which is representable by the form $x^2+hy^2$ can be factored into $p_1^2\cdots p_k^2q_1\cdots q_l$ where $p_i,q_j$ are irreducible elements in $R$ and $q_j$ is representable by $x^2+hy^2$ for all $j$.
\item[(3)] Some associate of a non-null prime element $p\in R$ is representable by $x^2+hy^2$ iff $-h$ is a square in $\mathbb{F}_{R/Rp}$, where $\mathbb{F}_{R/Rp}$  denotes the field of fractions of the quotient ring $R/Rp$.
\end{itemize}
\end{theorem}

The rest of the paper is structured as follows. In Section \ref{sec:ModCont} we define a generalization of the notion of continuant and describe some of its properties. Section \ref{sec:Q(x,y)ToQ(z,1)} is devoted to studying proper representations $x^2+gxy+h y^2$
in some commutative Euclidean rings, mainly in the ring of polynomials over a field of odd characteristic. In Section \ref{quadform} we consider proper representations $x^2+gxy+h y^2$ in the ring of integers. Some final remarks are presented in Section \ref{sec:final}.

\section{Generalized continuants}
\label{sec:ModCont}

With the aim of considering the problem of properly representing an element $m$ as
$x^2+gxy+h  y^2$, we extend the notion of continuants. Generalizations of continuants have previously appeared in the literature, mainly in commutative rings, where these generalizations can be considered determinants of certain matrices; see \cite[Sec.~8]{Urs58} and \cite{She72}. 

\begin{definition}[Generalized Continuants in Arbitrary Rings]\label{def:gCont}
In a ring $R$ we associate with the  element $[Q;h,s]$ the 3-tuple
formed from a sequence $Q$ of elements $q_1,q_2,\ldots,q_n$ of $R$, an element $h$ of $ R$
and an integer $s\ge 1$ via the following recurrence formula
\[
[q_1,\ldots q_n;h,s]=
\begin{cases}
[q_1,\ldots, q_n], &\text{if $s\ge n$;}\\
[q_1,\ldots, q_{n-1}]q_n +[q_1,\ldots, q_{n-2}]h,&\text{if $s=n-1$;}\\
[q_1,\ldots, q_{n-1};h,s]q_n+[q_1,\ldots, q_{n-2};h,s], & \text{if $s<n-1$.}
\end{cases}
\]
\end{definition}

This definition of generalized continuants carries several consequences, all of which are proved in Appendix \ref{sec:proofsCont}. These properties are referred to as Generalized Continuant Properties.
\begin{mylist}
\item\label{cuttingbis}
$[q_1,\ldots, q_n;h,s]=[q_1,\ldots, q_{s-1}]h[q_{s+2},\ldots, q_n]+
[q_1,\ldots, q_{s}][q_{s+1},\ldots, q_n]$, for $s<n$.
\item \label{zigzagbis}
If in a ring $R$ we find a unit $u$ commuting with each element $q_i$, then
\[
[u ^{-1}q_1,uq_2,\ldots, u  ^{(-1)^n}q_n;h,s]=\begin{cases}
[q_1,\ldots,q_n;h,s], &\text{for even $n$;}\\
u^{-1}[q_1,\ldots,q_n;h,s], &\text{for odd $n$.}
\end{cases}
\]
\end{mylist}

The next two properties pertain to commutative rings.
\begin{mylist}
\item\label{determibis}
 The generalized continuant $[q_1,\ldots, q_n;h,s]$
is the determinant of the tridiagonal $n\times n$ matrix $A=(a_{ij})$ with $a_{i,i}=q_i$
for $1\le i\le n$, $a_{i,i+1}=1$ for $1\le i<n$, $a_{s+1,s}=-h$ and $a_{i+1,i}=-1$ for $1\le i<n$
and $i\ne s$. See the determinant of the matrix below for a small example.
\[
[q_1,q_2,q_3,q_4,q_5;h,3]=\det\begin{bmatrix}
q_1& 1  &   &  &\\
  -1  &q_2&1  & & \\
    &-1  &q_3&1&\\
    & &-h&q_4&1\\
     &  &  &-1&q_5\end{bmatrix}
\]
\item\label{reverbis}
 $[q_1,q_2,\ldots,q_n;h,n-s]=[q_n,\ldots,q_2,q_1;h, s]$. \end{mylist}

\section{From   $Q(x,y)$ to $Q(z,1)$ and back}
\label{sec:Q(x,y)ToQ(z,1)}
 In this section, considering the form $Q(x,y)=x^2+gxy+hy^2$, we deal with the problem of going from a representation $Q(x,y)$ of an element $m$ to a multiple $Q(z,1)$ of $m$ and back. 
 
\subsection{From   $Q(x,y)$ to $Q(z,1)$} We begin with a general proposition which is valid for every commutative ring.

\begin{proposition}\label{direct}
If $Rx+Ry=R$ then there exist $z\in R$ such that $Q(x,y) $ divides $Q(z,1)$, where $Rm$ denotes the ideal generated by $m$.

If $R$ is Euclidean, we can explicitly find $z$  and the quotient $Q(z,1)/Q(x,y)$ with generalized continuants.
\end{proposition}

\begin{proof}
We have $u$ and $v$ such that $xu+yv=1$. Then,  computation with norms
in the ring obtained from $R$ by adjoining formally a root of the polynomial $T^2-gT+h$ provides the identity \[Q(x,y)Q(v-ug,u)=Q(xv-xug-yuh,xu+yv),\] which proves the first assertion. This identity can be interpreted also as a kind of Carroll's identity.

The determinant of the tridiagonal matrix
\begin{equation}\label{eq:det}
M=\begin{bmatrix}
q_s  & 1           &                 &           & & & &\\
-1     & \ddots  &\ddots &            & & &  &\\
         &\ddots   & q_2      &\ddots& & &  &\\
         &                 & -1         &q_1     &\ddots& &  &\\
        &                 &               & -h      &q_1+g &\ddots& &\\
        &                 &              &            &    -1    &q_2        &\ddots&\\
        &                 &             &              &          & \ddots & \ddots&1\\
       &                  &            &               &           &               &   -1 &    q_s
\end{bmatrix}
 \end{equation}
is $Q(x,y)$ with $x=[q_1,\ldots, q_s]$ and $y=[q_2,\ldots, q_s]$ if  $s\ge1$.

Also,
$Q(x,y)Q([q_1,\ldots q_{s-1}],[q_2,\ldots q_{s-1}])= Q(z,1)$,
where $z=(-1)^{s+1}c$ and $c$ is the determinant of the matrix formed by the $2s-1$ first rows and columns of $M$.
\end{proof}
The proof of Proposition \ref{direct} can be readily converted into a deterministic algorithm which finds a solution $z_0$ of $Q(z,1)\equiv0\pmod{m}$, given a representation $uQ(x,y)$ of an element $m$ in a Euclidean ring $R$ with a computable function $\nm$. See Algorithm \ref{DivAlg_FromQ(x,y)ToQ(z,1)}.
\begin{algorithm}[ht]
\SetKwData{Left}{left}\SetKwData{This}{this}\SetKwData{Up}{up}
\SetKwFunction{Union}{Union}\SetKwFunction{FindCompress}{FindCompress}
\SetKwInOut{Input}{input}\SetKwInOut{Output}{output}
\SetKw{Init}{initialisation:}
\SetKw{Ass}{assumptions:}
\Input{A commutative Euclidean ring $R$ with a computable function $\nm$.\\An element $m\in R$.\\A proper representation $uQ(x,y)$ of $m$, where $Q(x,y)=x^2+gxy+hy^2$.}
\Output{A solution $z_0$ of $Q(z,1)\equiv0\pmod{m}$ with $\nm(1)\le \nm(z_0)$.}
\tcc{Apply the Euclidean algorithm to $x$ and $y$ and obtain a sequence $(q_1,\ldots,q_s)$ of quotients.}
$s\leftarrow 0$\;
$m_0 \leftarrow m$\;
$r_0 \leftarrow z$\;
\Repeat{$r_{s}= 0$}
{
$s\leftarrow s+1$\;
$m_{s}\leftarrow r_{s-1}$\;
find $q_{s}, r_{s}\in R$ such that $m_{s-1} = q_{s} m_{s} +r_{s}$ with
$\nm(r_{s})<\nm(m_{s})$\;
}
$z_0 \leftarrow (-1)^{s+1}[q_s,q_{s-1},\ldots,q_1,q_1+g,q_2,\ldots,q_{s-1};h,s]$\;
\Return $z_0$
\caption{Deterministic algorithm for constructing a solution $z_0$ of $Q(z,1)\equiv0\pmod{m}$, given a proper representation $uQ(x,y)$ of an element $m$.}
\label{DivAlg_FromQ(x,y)ToQ(z,1)}
\end{algorithm}
\subsection{From $Q(z,1)$ to $Q(x,y)$} \label{subsec:Q(z,1)ToQ(x,y)}We begin the subsection with the following remark. 
\begin{remark}\label{rem:fChar2}
Let $R$ be a commutative ring.

If 2 is invertible, the form
$x ^2+gxy+hy^2$ can be rewritten as $(x+gy/2)^2+(h-g^2/4)y^2$. We may then assume $g=0$
without loss of generality.

If moreover $-h$ is an invertible square, say $h+k^2=0$, then
$x=\left(\frac{x+1}2\right)^2+h\left(\frac{x-1}{2k}\right)^2$
\end{remark}

Below we provide a proposition which can be considered as an extension of \cite[Prop.~16]{DP11}. 

\begin{proposition}\label{inverse}
Let $R=\F[X]$ be the ring of polynomials over  a field  $\F$  with characteristic different from 2,
and let $-h$ be a (non-null) non-square of $\F$.

If $m$ divides $z^2+ht^2$ with $z,t$ coprime, then $m$ is an associate of some $x^2+hy^2$ with $x,y$ coprime.
\end{proposition}
\begin{proof}
We introduce the extension $\G$ of $\F$
by a square root $\omega$ of $-h$. The ring $\G[X]$ is principal and $z^2+ht^2$ factorises as
$(z-\omega t)(z+\omega t)$. The two factors are coprime, since $2$ and $\omega$ are units, and any common divisor must divide their sum $2z$ and their difference $2\omega t$. 

Introduce
$x+\omega y=\gcd(m,z+\omega t)$. Then $x-\omega y$ is a gcd of $m$ and $z-\omega t$, using the
natural automorphism of $\G$. The polynomials $x-\omega y$ and $x+\omega y$ are coprime and both divide $m$. Thus, $m$ is divisible by $(x-\omega y)(x+\omega y)=x^2+hy^2$. On the other hand,  $m$ divides
$(z-\omega t)(z+\omega t)$. Consequently,  $m$ is an associate of  $x^2+hy^2$.
Since  $x-\omega y$ and $x+\omega y$ are coprime, $x$ and $y$ are also coprime.
\end{proof}

In the case of $m$ being prime, Proposition \ref{inverse} is embedded in Theorem 2.5 of the aforementioned paper of Choi et al. \cite{CLR1980}. 

\begin{proposition}\label{inversebis}
Let $R=\F[X]$ be the ring of polynomials over a field $\F$ with characteristic different from 2, and let  $h$ be a polynomial of degree 1.

If $m$ divides $z^2+ht^2$ with $z,t$ coprime, then $m$ is an associate of some $x^2+hy^2$ with $x,y$ coprime.
\end{proposition}
\begin{proof}
Consider the extension of the ring $R=\F[X]$ by a root of $T^2+h$; this extension of $R$ is isomorphic to $\F[T]$.

First assume that $h$ does not divide $m$. If $h$ and $z$ are not coprime, then $z$ can be rewritten as $z=z_1h$. Thus, $m$ dividing $h(hz_1^2+t^2)$ implies that $m$ divides $hz_1^2+t^2$, with $h$ and $t$ being coprime. Thus, we can assume that $h$ and $z$ are coprime.  Consequently,
$z^2+ht^2$ factors as $(z-Tt)(z+Tt)$, with the two factors being coprime.
 Reasoning as in Proposition~\ref{inverse}, we let $x+Ty$ be the gcd of $m$ and $z+Tt$ and we obtain that $x-Ty$ is the gcd of $m$ and $z-Tt$  and  that
$m$ is an associate of $x^2+hy^2$ with $x,y$ coprime.

If $m$ is a multiple of $h$, then $h$ does not divide $m/h$, since $z$ and $t$ are coprime. From the previous case it then follows that $m/h$  is an associate of some $x^2+hy^2$, with $x$ and $y$ coprime. Thus,  $m$ is an associate of $(hy)^2+h x^2$, with $hy$ and $x$ coprime.
\end{proof}

The next remark generalises \cite[Rem.~19]{DP11}.

\begin{remark}[Algorithmic considerations]\label{practical} For the cases covered in Pro\-po\-sitions \ref{inverse} and \ref{inversebis}, given an element $m$ and a solution $z_0$ of $z^2+h\equiv0\pmod{m}$, we can obtain a representation $x^2+hy^2$ of an associate of $m$ via generalized continuants and Brillhart's \cite{Bri72} optimisation. Indeed, divide $m$ by $z_0$ and stop when a remainder $r_{s-1}$ with degree at most $\deg(m)/2$ is encountered. This will be the $(s-1)$-th remainder, and $(uq_s,u^{-1}q_{s-1},\ldots,u^{(-1)^{s-2}}q_2)$ will be the quotients so far obtained. Then
\begin{align*}
x=\begin{cases}
r_{s-1},&\textrm{for odd $s$;}\\
u^{-1}r_{s-1},&\textrm{for even $s$.}
\end{cases}
\end{align*}
\begin{align*}
y=\begin{cases}
[uq_s,u^{-1}q_{s-1},\ldots,u^{(-1)^{s-2}}q_2],&\textrm{for odd $s$;}\\
u^{-1}[uq_s,u^{-1}q_{s-1},\ldots,u^{(-1)^{s-2}}q_2],&\textrm{for even $s$.}
\end{cases}
\end{align*}
This conclusion follows from dividing 
\begin{align*}
m/u&=[q_s,\ldots,q_1,q_1,\ldots,q_s;h,s]\quad \text{by}\\
z_0&=[q_{s-1},\ldots,q_1,q_1,\ldots,q_s;h,s-1],
\end{align*}
using generalized continuant properties.
\end{remark}
Remark \ref{practical} can be readily translated into a deterministic algorithm for computing representations $Q(x,y)$;  see Algorithm \ref{PolyDivAlg}. 

The argument presented in \cite[Prop.~17]{DP11} can be applied to the form $x^2+h y^2$ in polynomials over a  field $\F$ of characteristic different from 2, where $-h$ is either a non-square $\in \F$ or a polynomial in $\F[X]$ of degree 1. This argument implicitly invokes the uniqueness of the quotients and the remainders in the division algorithm.

\begin{corollary}[of Proposition \ref{inverse}: {$-h$ a non-square unit in $\F[X]$}]\label{applinvbis}
 Let $m$ be a non-unit of $\F[X]$ and a divisor of $z^2+h$ for some $z\in \F[X]$ with $\deg(z)<\deg(m)$. Then, $m=(x^2+h y^2)u$ for some unit $u$ and the Euclidean algorithm on $m$ and $z$ gives the unit $u$ and the sequence  \[(uq_s,u^{-1}q_{s-1},\ldots,u^{(-1)^{s+1}}q_1,u^{(-1)^{s}}h^{-1}q_1,\ldots, u^{-1}h^{(-1)^s}q_s)\] such that $x=[q_1,\ldots q_s]$ and $y=[q_2,\ldots,q_s]$.
\end{corollary}

In the next example we illustrate Remark \ref{practical} and the method of Corollary \ref{applinvbis}, in this order, for the case of $\F=\mathbb{Q}$. Let $h=3$ and $m=1+2X+3X^2+2X^3+X^4$. Then $m$ divides  $((5+12X+6X^2+4X^3)/3)^2+3$. The Euclidean division gives
\begin{align*}
1+2X+3X^2+2X^3+X^4=&((5+12X+6X^2+4X^3)/3)(3X/4+3/8)\\
&{}+3/8-3X/4-3X^2/4.
\end{align*}
Here the first remainder has degree at most $\deg(m)/2$, thus we stop the division process and obtain
$s=2$, $x=(3/8-3X/4-3X^2/4)/u$ and $y=[3X/4+3/8]/u$. It is now routine  to get $u=9/16$.

If instead we use the method of Corollary \ref{applinvbis}, then we obtain the unit $u=9/16$ and the sequence \[(9/16\cdot2/3\cdot (1 + 2 X), 16/9\cdot (-1/2 - X), 9/16\cdot1/3\cdot (-1/2 -X),  16/9\cdot3\cdot2/3\cdot (1 + 2 X)).\] From this sequence we conclude that
\begin{align*}
x&=[2/3\cdot (1 + 2 X), -1/2 - X],\\
y&=[2/3\cdot (1 + 2 X)]. 
\end{align*}
\begin{corollary}[of Proposition \ref{inversebis}: {$h$ of degree 1 in $\F[X]$}]  \label{applinvter}
Let $m$ be a po\-ly\-nomial over $\F[X]$ and a divisor of $z^2+h$ for some $z\in \F[X]$ with $\deg(z)<\deg(m)$ and $z,h$ coprime. Then, $m=(x^2+h y^2)u$ for some unit $u$ and the values of $x$ and $y$ can be obtained by Remark \ref{practical}.
\end{corollary}

Consider the following example. Let $h=X$, $m=1 +  X + X^3 + X^4$ and $z=(X^3+2X^2+1)/2$. Then, the division gives
\begin{align*}
1 +  X + X^3 + X^4=&(X^3+2X^2+1)/2\cdot(-2+ 2 X)+2+2X^2.
\end{align*}

At this step we should stop the division process as the first remainder has
at most half the degree of $m$. Now we know that $s=2$, $x=(2 + 2 X^2)/u$ and $y=u^{-1}[-2 + 2 X]$ for a unit $u$. It plainly follows that $u=4$.

\begin{algorithm}[ht]
\SetKwData{Left}{left}\SetKwData{This}{this}\SetKwData{Up}{up}
\SetKwInOut{Input}{input}\SetKwInOut{Output}{output}
\SetKw{Init}{initialisation:}
\SetKw{Ass}{assumptions:}
\Input{A field $\F$ with characteristic different from 2.\\The ring $R=\F[X]$  of polynomials over $\F$.\\A square-free element $h\in \F$ or a polynomial $h\in R$ of degree 1.\\A polynomial  $m$ with $\nm(1)< \nm(m)$.\\A solution $z_0$ of $Q(z,1)\equiv0\pmod{m}$ with $\nm(1)< \nm(z_0)<\nm(m_0)$.}
\Output{A unit $u$ and a proper representation $uQ(x,y)$ of $m$.}
\Ass{} The polynomials $z$ and $h$ are coprime.

\tcc{Divide $m$ by $z$ using the Euclidean algorithm until we find a remainder $r_{s-1}$  with degree at most $\deg(m)/2$.}
$s\leftarrow1$\;
$m_0\leftarrow m$\;
$r_0\leftarrow z$\;
\Repeat{$\deg(r_{s-1})\le \deg(m)/2$}
{
$s\leftarrow s+1$\;
$m_{s-1}\leftarrow r_{s-2}$\;
find $k_{s-1}, r_{s-1}\in R$ such that $m_{s-2} = k_{s-1} m_{s-1} +r_{s-1}$ with
$\nm(r_{s-1})<\nm(m_{s-1})$\;
}
\tcc{Here we have a sequence $(k_1,\ldots,k_{s-1})$ of quotients.}
$x_{temp}\leftarrow r_{s-1}$\;
$y_{temp}\leftarrow [k_{1},\ldots,k_{s-1}]$\;
\tcc{We obtain a unit $u$.}
\lIf{$s$ is odd} {Solve $m=u(x_{temp}^2+hy_{temp}^2)$ for $u$} \\ \lElse{Solve $um=x_{temp}^2+hy_{temp}^2$ for $u$}

\tcc{We obtain $(x,y)$ so that  $m=(x^2+hy^2)u$.}
\lIf{$s$ is odd}{$x\leftarrow x_{temp}$} \lElse{$x\leftarrow u^{-1}x_{temp}$}\;
\lIf{$s$ is odd}{$y\leftarrow y_{temp}$} \lElse{$y\leftarrow u^{-1}y_{temp}$}\;
\Return $(x,y,u)$
\caption{Deterministic algorithm for constructing a proper representation $uQ(x,y)=u(x^2+hy^2)$ of an element $m$.}
\label{PolyDivAlg}
\end{algorithm}

We  may now wonder how far can we push this method for polynomials over a field of characteristic different from 2? That is, will the method work for $h$ with $\deg(h)>1$ over any such field?

 We first note that the property

``$m | z^2+h \Rightarrow \exists x,y,u \;(m= u(x^2+y^2h) \land u \text{~unit})$''

\noindent does not hold in general for  $h$ reducible. Consider $h=X^3+X^2+X$ in polynomials over a field of characteristic $\ne 3$. Then, $X^2+X+1$ divides $0^2+1^2h$ and is certainly not of the form $x^2+y^2h$. Indeed, here we have either $y^2h$ null or of odd degree $\ge 3$. In the former case,
it follows that $x^2+y^2h =x^2$ is  a square, but $X^2+X+1$ is not
       a square in a field of characteristic $\ne 3$, while in the latter case $x^2+y^2h$ has degree
$\ge 3>\text{degree}(X^2+X+1)$.

What about irreducible $h$ with $\deg(h)\ge2$? Already for degree 2
the property does not hold in general. Indeed, consider in $\Q[X]$ the polynomials $h=X^2-2$, $z=X^2$ and $m=X-1$. Observe that $X-1$ divides $X^4+X^2-2$ and that $X-1$ is not of the form $u(x^2+y^2(X^2-2))$. To see this note that the degree of $x ^2+y^2(X^2-2)$ is either $2\deg(x)$ or $2+2\deg(y)$.

For specific fields we find situations where the property holds. Take, for instance, the field $\mathbb{R}$ of reals, $h=X^2+1$ and every real polynomial $m$ taking only positive values over $\mathbb{R}$. It is known that any polynomial $m$ over $\mathbb{R}$, which takes at every point of $\mathbb{R}$ a  positive value, has the form $\prod(a_k X^2+2b_kX+c_k)$, where $a_k,b_k,c_k\in \mathbb{R}$ and $b_k^2-a_kc_k<0$. Thus, it suffices to consider the case of $m=aX^2+2bX+c$ with $a>0,c>0$ and $b^2-ac<0$.
If $b=0$, then
\[m=\begin{cases}(\sqrt{a-c}X)^2+\sqrt c ^2(X^2+1),& \textrm{if $a\ge c$;}\\ \sqrt{c-a}^2+\sqrt a ^2(X^2+1),& \textrm{if $a\le c$.}  \end{cases}\]

If instead $b\ne0$, then, setting $d=\sqrt{(a+c-\sqrt{(a-c)^2+4b^2})/2}$, we obtain $m=(\sqrt{a-d^2}X+e \sqrt{c-d^2})^2+d^2(X^2+1)$, where $e=\pm 1$ is 
the sign of $b$.

For $h=-X^2-1$ and every real polynomial $m$ over $\mathbb{R}$, we have another situation where the property holds. Observe that the form $Q(x,y)=x'^2+(-X^2-1)y'^2$ is equivalent to the form $Q(x,y)=x^2+2Xxy-y^2$ (by Remark \ref{rem:fChar2}). Any polynomial of degree 1 is an associate of some $a^2-b^2+2abX$ with units $a$ and $b$. We now take care of polynomials $m=k(X^2+2vX+w)$ with no real zeros and $k,v,w\in \mathbb{R}$. Here note that $v^2<w$. Set  $p(X)=(X+a)^2+2b(X+a)X-b^2$. We solve the equation $p(X)=m$ in $(a,b,k)$. We first find that
$k=1+2b$, $a=kv/(1+b)$ (if $b\ne-1$) and $w=(a^2-b^2)/k=-X^2-2vX$. If $b=-1$ then $v=0$, $k=-1$ and $a=\pm\sqrt{1-w}$ with $0<w\le1$. If instead $b\ne-1$, then, substituting $a=kv/(1+b)$ into $-(1+b)^2p(X)$, we obtain \[b^4+2(w+1)b^3+(5w-4v^2+1)b^2+4(w-v^2)b+w-v^2=0.\]

This expression in $b$ is $1/16$ when $b=-1/2$ and $-v^2$ when $b=-1$. Hence there is a solution $b$ in the open interval $(-1,-1/2)$ for $w>1$. Consequently, each real polynomial is an associate of some polynomial $x^2+2xyX-y^2=(x+Xy)^2+(-1-X^2)y^2$.

Using the automorphisms of $ \mathbb{R}[x]$, both previous approaches can easily be applied to any real polynomial of degree 2 with no real roots.

\section{From $Q(z,1)$ to $Q(x,y)$: integral quadratic forms}\label{quadform}

In this section, given  integers $m,z$ such $m | Q(z,1)$, we provide an algorithm that proves the existence of representations $Q(x,y)u$ of $m$ for a unit $u$ and certain forms $Q$.

Since 2 is not invertible in $\Z$, we have to consider the rings of algebraic integers of $\Q[\sqrt{-h}]$, that is, the rings $\Z[\sqrt{-h}]$ for forms $x^2+h y^2$ with $|h|$ square-free and $h\not\equiv -1\pmod 4$, and the ring $\Z[(1+\sqrt{1-4h})/2]$ for forms $x^2+xy+ h y^2$ with $|1-4h|$ square-free; see \cite[p.~35]{Sam70}.

What are those rings of integers for which the following property holds?

``$m | Q(z,1) \Rightarrow \exists x,y,u \;(m= uQ(x,y) \land u \text{~unit})$.''

The answer is given by the rings whose corresponding forms have  class number $H(\Delta)$ equal to one \cite[pp.~6--7, pp.~81--84]{bue89}. Here $\Delta$ denotes the form discriminant. In the case of $\Delta<0$, all the principal rings satisfy the property; these
values of $\Delta$ are the following: $-3,-4,-7,-8,-11,-19,-43,-67,-163$; see \cite[A014602]{OEIS} and \cite[pp.~81]{bue89}. For these negative fundamental discriminants, generalized continuants provide a constructive proof of the property. It is also known that there are four negative non-fundamental discriminants of class number one, namely $-12, -16, -27$ and $-28$; see \cite[pp.~81]{bue89} and \cite[Thm.~7.30]{Cox89}. For the case of $\Delta>0$, while we do not even know whether the list of such determinants is infinite, it is conjectured this is the case; see \cite[pp.~81--82]{bue89} and \cite[A003655]{OEIS}.

Recall the class number $H(\Delta)$ \cite[p.~7]{bue89} gives  the number of equivalence classes of integral binary quadratic forms with discriminant $\Delta \in \mathbb{Z}$.

Next we recall the following well-known result of Rabinowitsch \cite{Rab1913}.
\begin{theorem}[{\cite{Rab1913}}]\label{thm:Rabi} For a fundamental discriminant $\Delta=1-4\kappa\le-7$, it follows that $H(\Delta)=1$ iff $x^2+x+\kappa$ attains only prime values for $-(\kappa-1)\le x\le\kappa-2$. 

\end{theorem}
Below  we present a division algorithm (Algorithm \ref{DivAlgI}) which, for any negative fundamental discriminant of class number one, gives a proper representation $Q(x,y)$ of $m$, provided that $m$ divides $Q(z,1)$.  

Since $m=1$ trivially admits the proper representation $(1,0)$ of $Q(x,y)=x^2+gxy+hy^2$,
Algorithm \ref{DivAlgI} assumes $|m|>1$.
\begin{algorithm}[ht]
\SetKwData{Left}{left}\SetKwData{This}{this}\SetKwData{Up}{up}
\SetKwFunction{Union}{Union}\SetKwFunction{FindCompress}{FindCompress}
\SetKwInOut{Input}{input}\SetKwInOut{Output}{output}
\SetKw{Init}{initialisation:}
\SetKw{Ass}{assumptions:}
\LinesNumbered
\Input{A negative fundamental discriminant $\Delta$ of class number one.\\An integer $m_0$ with $1< |m_0|$.\\ A solution $z_0$ of $Q(z,1)\equiv0\pmod{m_0}$ with $1< |z_0|<|m_0|$.}
\Output{A proper representation $Q(x,y)$ of $m_0/u$ ($u=\pm1$).}
$s\leftarrow0$\;
\While{$|m_{s}|\ne1$}
{
$s\leftarrow s+1$\;
$m_s\leftarrow Q(z_{s-1},1)/m_{s-1}$\;
find $k_s, z_s\in R$ such that $z_{s-1} = k_sm_s +z_s$ with
$|z_s|<|m_s|$\;
\tcc{We prioritise  non-null quotients $k_s$.}
}
\tcc{Here we have the unit $m_s$ and sequence $(k_1,\ldots,k_s)$.}
\tcc{To keep consistency with the previous sections  of the paper we reverse the subscripts of the quotients}
$(q_1,q_{2},\ldots,q_s)\leftarrow(k_s,k_{s-1},\ldots,k_1)$\;
$x\leftarrow[m_sq_1,m_s^{-1}q_{2},\ldots,m_s^{(-1)^{s-2}}q_{s-1},m_s^{(-1)^{s-1}}q_s]$\;
$y\leftarrow[m_s^{-1}q_{2},\ldots,m_s^{(-1)^{s-2}}q_{s-1},m_s^{(-1)^{s-1}}q_s]$\;
\Return $(x,y)$
\caption{Deterministic algorithm for constructing a proper representation $Q(x,y)=x^2+gxy+hy^2$ of an element $m$.}
\label{DivAlgI}
\end{algorithm}
\begin{remark}[Algorithm \ref{DivAlgI}: Prioritising non-null quotients]
In the Euclidean division of $z_{s-1}$ by $m_s$ with $|z_{s-1}|<|m_s|$, a valid quotient $k_s$ could be $\pm1$ or 0. By ``prioritising non-null quotients $k_s$'' we mean that, in this situation, we always choose the  $k_s$ which is not null.
\end{remark}
\begin{proposition}[Algorithm \ref{DivAlgI} Correctness]\label{correctness} Let $h,g,\Delta,u,m_0,z_0$, and  $m_i,z_i$, $q_i$ ($i=1,\ldots,s$) be as in Algorithm \ref{DivAlgI}. Then, Algorithm \ref{DivAlgI} produces a proper representation $Q(x,y)$ of $m_0/u$.
\end{proposition}
\begin{proof} In the ring $\Z[(1+\sqrt{1-4h})/2]$ we consider the form $Q(x,y)=x^2+xy+hy^2$, while in the ring $\Z[\sqrt{-h}]$ we consider the form $Q(x,y)=x^2+hy^2$. As the proof method is the same in both cases, we restrict ourselves to the former case, that is, to the case of $\Delta=-3,-7,-11,-19,-43,-67,-163$ and $h=1, 2, 3, 5, 11, 17, 41$.

\par\noindent{\bf Claim 1. Algorithm \ref{DivAlgI} terminates with the last $m_j$ being $\pm1$.}

 As a general approach we show that the sequence $|m_i|$ ($i=0,\ldots, s-1$) is decreasing, that is, $|m_{i+1}|<|m_i|$ . Once this decreasing character fails, we show that the algorithm  stops with the last $m_i$ being a unit.

Recall we have $|m_i|\ge|z_i|+1$ for $i=1,\ldots, s-1$.

 {\bf Case $(\Delta,h)=(-3,1),(-7,2)$:} $|m_i||m_i|\ge z_i^2+2|z_i|+1> |z^2_i+z_i+h|=|m_i||m_{i+1}|$, and thus $|m_i|>|m_{i+1}|$ for $|z_i|>1$. Assume that for a certain $z_i$, say $z_{s-1}$, $|z_{s-1}|=1$. If $h=1$ then Line 4 ($z_{s-1}^2+z_{s-1}+1=m_{s-1}m_{s}$) gives that $|m_s|=1$, as desired. In the case of $h=2$ and $z_{s-1}=-1$, we have that $z_{s-1}^2+z_{s-1}+2=2$ and $|m_s|=1$. If $h=2$ and $z_{s-1}=1$,  we have from Line 4 again that $1+1+2=m_{s-1}m_s$. Then we deduce that either $|m_{s-1}|=2$ and $|m_{s}|=2$ or $|m_{s-1}|=4$ and $|m_{s}|=1$.  The configuration $|m_{s-1}|=4$ and $|m_{s}|=1$ will cause the algorithm to stop with $m_s$ being a unit. In the case of $|m_{s-1}|=2$ and $|m_{s}|=2$, in Line 5 we have $z_{s-1}=k_{s}m_{s}+z_{s}$ and the algorithm would give $z_{s}=-1$, which implies $|m_{s+1}|=1$.

Consequently, in these two cases Algorithm \ref{DivAlgI} terminates with the last $m_j$ being a unit.

{\bf Case $(\Delta,h)=(-11,3),(-19,5),(-43,11),(-67,17),(-163,41)$:} Reasoning as in the previous case, we have that $|m_i||m_i|\ge z_i^2+2|z_i|+1>|z^2_i+z_i+h|=|m_i||m_{i+1}|$, unless $|z_i|\le h-1$. Suppose $|z_{s-1}|\le h-1$. Note that $|m_{s-1}|>1$, otherwise the algorithm would have stopped. By Theorem \ref{thm:Rabi}, the polynomial $Q(z_{s-1},1)=z_{s-1}^2+z_{s-1}+h$ is prime for $-(h-1)\le z_{s-1}\le h-2$. Thus, we have that $|m_{s-1}|>|m_s|$ with $|m_s|=1$.

If instead $z_{s-1}= h-1$, then $|m_{s-1}|=|m_{s}|=h$, which implies that $z_{s}=-1$ and $|m_{s+1}|=1$, causing the algorithm to stop.

\par\noindent{\bf Claim 2. Algorithm \ref{DivAlgI} produces  a proper representation $Q(x,y)$.}

It only remains to prove that $x$ and $y$ have the required form. First we reverse the subscripts
of the quotients, that is, the quotient $k_s$ becomes $q_1$, the quotient $k_{s-1}$ becomes $q_2$, and so on. Thus, after the ``while'' loop we have that $z_{s-1}=m_sq_1$, where $m_s$ is a unit. We know that $Q(m_sq_1,1)=m_{s-1}m_s$. Consequently,
$m_{s-1}m_s=[m_sq_1+1,m_sq_1;h,1]$. Then, by Property P--\ref{reverbis} and Property P--\ref{cuttingbis} it follows
$$z_{s-2}=[m_s^{-1}q_{2},m_sq_{1}+1,m_sq_1;h,2].$$
Then, from the equation $m_{s-2}m_{s-1}=Q(z_{s-2},1)$ we obtain
\[m_{s-2}m_s^{-1}=[m_s^{-1}q_{2},m_sq_{1}+1,m_sq_1,m_s^{-1}q_{2};h,2]=Q([m_sq_1,m_s^{-1}q_2],[m_s^{-1}q_2]).\]

Continuing this process, we have
\begin{align*}
z_{0}&=[m_s^{(-1)^{s-1}}q_s,K,m_sq_1+1,m_sq_1,K^{-1};h,s],\\
m_0m_s^{(-1)^{s+1}}&=[m_s^{(-1)^{s-1}}q_s,K,m_sq_1+1,m_sq_1,K^{-1},m_s^{(-1)^{s-1}}q_s;h,s],
\end{align*}
where $K=m_s^{(-1)^{s-2}}q_{s-1},\ldots,m_s^{-1}q_2$ and $K^{-1}=m_s^{-1}q_2,\ldots,m_s^{(-1)^{s-2}}q_{s-1}$.

Consequently, from Property P--\ref{zigzagbis} it follows that $m_0m_s^{(-1)^{s+1}}=Q(x,y)$, where
\begin{align*}
x=&[m_sq_1,m_s^{-1}q_{2},\ldots,m_s^{(-1)^{s-2}}q_{s-1},m_s^{(-1)^{s-1}}q_s],\\
y=&[m_s^{-1}q_{2},\ldots,m_s^{(-1)^{s-2}}q_{s-1},m_s^{(-1)^{s-1}}q_s].
\end{align*}
Using Property P--\ref{reverbis} we can write $m_0m_s^{(-1)^{s+1}}$ as follows; see Equation (\ref{eq:det}). 
\begin{align*}
m_0m_s^{(-1)^{s+1}}&=[m_s^{(-1)^{s-1}}q_s,K,m_sq_1,m_sq_1+1,K^{-1},m_s^{(-1)^{s-1}}q_s;h,s].
\end{align*}
\begin{remark} \label{class1} Note that we require all the numbers $m_i$ to be represented by the form $Q(x',y')$. This is assured by the fact that $Q$ has class number one.
\end{remark}
As a result, each root $z_0$ of $Q(z,1)\equiv0\pmod{m_0}$ with $ 1< |z_0|< |m_0|$ gives rise to a proper representation of $m_0m_s^{(-1)^{s+1}}$ as $Q(x,y)$. The coprimality of $x$ and $y$ follows from Property P--\ref{lewiscarroll}.
\end{proof}
Let us see an example. For the ring $\Z[(1+\sqrt{-19})/2]$ the form is $x^2+xy+5y^2$. Take $m_0=251$ and $z_0=52$. Then $251\cdot 11=52^2+52+5$, and the division gives
\begin{align*}
251\cdot 11&=52^2+52+5 &\rightarrow& & 52&=4\cdot 11 +8,\\
11\cdot 7&=8^2+8+5 &\rightarrow& &  8&=1\cdot 7 +1,\\
7\cdot 1 &= 1^2+1+5 &\rightarrow& & 1&=1\cdot 1.
\end{align*}
Thus, we have $m_3=1$ and $(q_3,q_2,q_1)=(4,1,1)$. From this we
recover the continuant representation of $m_0=251$
\[
251=\det\begin{bmatrix}
4&1&&&&\\-1&1&1&&&\\&-1&1&1&&\\&&-5&1+1&1&\\&&&-1&1&1\\&&&&-1&4
\end{bmatrix}.
\]
Consequently, we conclude that $251=x^2+xy+5y^2$ with $x=[1,1,4]=9$ and $y=[1,4]=5$.

\section{Final remarks}
\label{sec:final}

In Algorithm \ref{DivAlgI} we require $m_s$ to be $\pm1$. However, this may be an unnecessarily strong restriction.
 If in Algorithm \ref{DivAlgI} we replace the condition of the while loop by $z_s\ne0$, then this modified Algorithm \ref{DivAlgI} may also
 end with the last $m_j$, say $m_s$, being different from $\pm1$. Further, if such  $m_s$ admits a representation as $Q(x,y)$,
 then the formula
\begin{align}
(x^2 + gxy + hy^2)(z^2 + gzw +hw^2)=&{}(xz-hyw)^2+g(xz-hyw)\times\label{eq:Sq}\\
&{}\times(xw+yz+gyw)+\nonumber\\
&{}+h(xw+yz+gyw)^2\nonumber
\end{align}
will provide a desired representation of $m=m_0$ for a larger number of forms $Q(x,y)$.
First recall that in Algorithm \ref{DivAlgI}
\begin{align*}
m_0m_s^{(-1)^{s+1}}=&{}x^2+gxy+hy^2,
\end{align*}
where \begin{align*}
x=&[m_sq_1,m_s^{-1}q_{2},\ldots,m_s^{(-1)^{s-2}}q_{s-1},m_s^{(-1)^{s-1}}q_s],\\
y=&[m_s^{-1}q_{2},\ldots,m_s^{(-1)^{s-2}}q_{s-1},m_s^{(-1)^{s-1}}q_s].
\end{align*}
Then, to recover the representation of $m_0$ (associated with $z_0$) we just need to express $m_s$ or $-m_s$ as $Q(x,y)$. This simple modification of Algorithm \ref{DivAlgI} will provide proper representations $Q(x,y)$ of $\pm m$ for some forms $Q$ with discriminant $\Delta>0$ and $H(\Delta)=1$; see \cite{Disc_Real_Quad_Field_NarrowClassNo1}. The following two examples illustrate this idea. Recall that the condition of the ``while'' loop is now $z_s\ne0$. 

For the ring $\Z[(1+\sqrt{17})/2]$ the form is $x^2+xy-4y^2$. Take $m_0=3064$ and $z_0=564$.  Noticing $3064\cdot 104=564^2+564-4$, the division gives
\begin{align*}
3064\cdot 104&=564^2+564-4 &\rightarrow& & 564&=5\cdot 104 +44,\\
104\cdot 19&=44^2+44-4 &\rightarrow& &  44&=2\cdot 19 +6,\\
19\cdot 2 &= 6^2+6-4 &\rightarrow& & 6&=3\cdot 2.
\end{align*}
Thus, we have $s=3$, $m_3=2$ and $(q_3,q_2,q_1)=(5,2,3)$. From this we
recover the continuant representation of $m_0\cdot 2=6128$
\[
6128=\det\begin{bmatrix}
10&1&&&&\\-1&1&1&&&\\&-1&6&1&&\\&&4&6+1&1&\\&&&-1&1&1\\&&&&-1&10
\end{bmatrix}.
\]

The  representation $Q(x,y)$ of 6128 is given by $x=[2\cdot3,2^{-1}\cdot2,2\cdot5]=76$ and $y=[2^{-1}\cdot2,2\cdot5]=11$. Note that 
$2=2^2+2\cdot1-4\cdot1^2$. Using Equation (\ref{eq:Sq}) in the form $3064\cdot (2^2+2\cdot 1-4\cdot 1)=76^2+76\cdot 11-4\cdot 11^2$, we conclude that $3064=92^2+92\cdot (-27)-4(-27)^2$.

For the ring $\mathbb{Z}[\sqrt{6}]$ the form is $Q(x,y)=x^2-6y^2$. Take $m_0=37410$ and $z_0=1326$. Noticing $37410\cdot 47=1326^2-6$, the division gives
\begin{align*}
37410\cdot 47&=1326^2-6 &\rightarrow& & 1326&=28\cdot (47) +10,\\
47\cdot 2&=10^2-6 &\rightarrow& &  10&=5\cdot 2.
\end{align*}
Thus, we have $s=2$, $m_2=2$ and $(q_2,q_1)=(28,5)$. The  representation $Q(x,y)$ of $37410\cdot 2^{-1}$ is $x=[2\cdot 5, 2^{-1}\cdot 28]=141$ and $y=[2^{-1}\cdot 28]=14$. Note that $-2=2^2-6\cdot1^2$. Using Equation (\ref{eq:Sq}) in the form $(141^2-6\cdot 14^2)(2^2-6\cdot1^2)=-37410$, we have that $-37410=366^2-6\times 169^2$.
 
Below we present some of the forms $Q(x,y)$ for which the proposed modification of Algorithm \ref{DivAlgI} will give the representation of $m_0$ associated with the given $z_0$.  The forms are given in the format  $(Q,\{\textrm{list of possible values of $m_s$}\})$. Note that $m_s$ is a divisor of $h$ and that, for every case, either $m_s$ or $-m_s$ is represented by the form.

$\begin{array}{ll}
(x^2-2y^2,\{\pm1,\pm2\})&(x^2+xy-4y^2,\{\pm1,\pm2,\pm4\})\\
(x^2-3y^2,\{\pm1,\pm3\}) &(x^2+xy-7y^2,\{\pm1,\pm7\}) \\
(x^2-6y^2,\{\pm1,\pm2,\pm3,\pm6\})&(x^2+xy-9y^2,\{\pm1,\pm3,\pm9\})\\
(x^2-7y^2,\{\pm1,\pm7\})&(x^2+xy-10y^2,\{\pm1,\pm2,\pm5,\pm10\})\\
(x^2+xy-y^2,\{\pm1\})&(x^2+xy-13y^2,\{\pm1,\pm13\})\\
(x^2+xy-3y^2,\{\pm1\})&(x^2+xy-15y^2,\{\pm1,\pm3,\pm5,\pm15\})
\end{array}$

Recall that Matthews \cite{Mat02} provided representations of certain integers as $x^2-hy^2$, where $h = 2, 3, 5, 6, 7$. From the previous remarks it follows that the modified Algorithm \ref{DivAlgI} covers
all the cases studied by Matthews \cite{Mat02}. Observe that the form $x^2-5y^2$, studied by Matthews \cite{Mat02} and associated with the non-principal ring $\Z[\sqrt{5}]$, has been superseded by the form $x^2+xy-y^2$ associated with the integral closure of $\Z[\sqrt{5}]$, that is, $\Z[(1+\sqrt{5})/2]$. It is known that the forms $x^2+xy-y^2$ and $x^2-5y^2$ represent the same integers. Indeed, if an integer $m$ is represented by $x^2-5y^2$ then the identity $x^2-5y^2=(x-y)^2+(x-y)(2y)-(2y)^2$ gives a representation, not necessarily proper, of $m$ by the form $x'^2+x'y'-y'^2$. If instead an integer $m$ is represented by the form $x^2+xy-y^2$, then, depending on the parity of $x$ and $y$, one of the identities 
\begin{align*}
x^2+xy-y^2&=\left(x+y+\frac{x+2y}{2}\right)^2-5\left(\frac{x+2y}{2}\right)^2,\\
&=\left(2x-y+\frac{-x+y}{2}\right)^2-5\left(\frac{-x+y}{2}\right)^2,\\
&=\left(x+\frac{y}{2}\right)^2-5\left(\frac{y}{2}\right)^2
\end{align*}
gives a representation by the form $x'^2-5y'^2$. 

Unsatisfactorily, our algorithm does not terminate for every $\Delta>0$ with $H(\Delta)=1$. For instance, take $\Delta=73$, $m_0=267$ and $z_0=23$. The corresponding quadratic form is $x^2+xy-18y^2$, and we have that $267=(-69)^2+(-69)\cdot14-18\cdot14^2$ and $267|23^2 + 23 - 18$.

The approach presented in the paper is likely to work for other representations if new generalized continuants are defined.

{\it Mathematica$^\circledR$}\cite{Wol99} implementations of most of the algorithms presented in the paper and other
related algorithms are available at

\url{http://guillermo.com.au/wiki/List_of_Publications}

\noindent under the name of this paper.
\section{Acknowledgments}

The authors thank the referee for his/her careful and thoughtful review. The paper's presentation has certainly benefited from his/her comments.  
\appendix 

\section{Proofs of some of the generalized continuant properties} 
\label{sec:proofsCont}
\begin{proposition}[Property P--\ref{cuttingbis}]\label{prop:P3} For integers $n,h,s$ such that $1\le s<n$ and elements $q_1,\ldots, q_n$ of a ring $R$, the following identity holds: \[[q_1,\ldots, q_n;h,s]=[q_1,\ldots, q_{s-1}]h[q_{s+2},\ldots, q_n]+
[q_1,\ldots, q_{s}][q_{s+1},\ldots, q_n].\]
\end{proposition}
\begin{proof} In the case of $h=1$ generalized continuants reduce to the traditional continuants and Property P--\ref{cuttingbis} reduces to the  well-known identity
\[[q_1,\ldots, q_n]=[q_1,\ldots, q_{s-1}][q_{s+2},\ldots, q_n]+
[q_1,\ldots, q_{s}][q_{s+1},\ldots, q_n].\]
This identity is proved for commutative rings in \cite[Lem.~1]{CELV99} and  \cite[Sec.~6.7]{GKP}, but the approach by Graham et al.~\cite[Sec.~6.7]{GKP} works for any ring.

For any $h>1$ and $s=n-1$, Property P--\ref{cuttingbis} follows from the definition of generalized continuants.

Consider any $h>1$ and $s=n-2$. Then, from Definition \ref{def:gCont} it follows that
\begin{align*}
[q_1,\ldots q_n;h,n-2]&=[q_1,\ldots, q_{n-1};h,n-2]q_n+[q_1,\ldots, q_{n-2};h,n-2],\\
&=([q_1,\ldots,q_{n-3}]h+[q_1,\ldots, q_{n-2}][q_{n-1}])q_n+[q_1,\ldots, q_{n-2}],\\
&=[q_1,\ldots,q_{n-3}]hq_n+[q_1,\ldots, q_{n-2}]([q_{n-1}])q_n+1),\\
&=[q_1,\ldots,q_{n-3}]hq_n+[q_1,\ldots, q_{n-2}][q_{n-1},q_n].
\end{align*}
Finally, fix $h>1$ and $s<n-2$ and proceed by induction on $n$. The base cases $n=1,2$ fall in the previous cases. From Definition \ref{def:gCont} it follows that
\[[q_1,\ldots q_n;h,s]=[q_1,\ldots, q_{n-1};h,s]q_n+[q_1,\ldots, q_{n-2};h,s].\]
By the induction hypothesis we have the following.
\begin{align*}
[q_1,\ldots, q_{n-1};h,s]&=[q_1,\ldots, q_{s-1}]h[q_{s+2},\ldots, q_{n-1}]+
[q_1,\ldots, q_{s}][q_{s+1},\ldots, q_{n-1}];\\
[q_1,\ldots, q_{n-2};h,s]&=[q_1,\ldots, q_{s-1}]h[q_{s+2},\ldots, q_{n-2}]+
[q_1,\ldots, q_{s}][q_{s+1},\ldots, q_{n-2}].
\end{align*}
Thus
\begin{align*}
[q_1,\ldots q_n;h,s]&=([q_1,\ldots, q_{s-1}]h[q_{s+2},\ldots, q_{n-1}]+
[q_1,\ldots, q_{s}][q_{s+1},\ldots, q_{n-1}])q_n\\
&\quad+[q_1,\ldots, q_{s-1}]h[q_{s+2},\ldots, q_{n-2}]+
[q_1,\ldots, q_{s}][q_{s+1},\ldots, q_{n-2}],\\
&=[q_1,\ldots, q_{s-1}]h\left([q_{s+2},\ldots, q_{n-1}]q_n+[q_{s+2},\ldots, q_{n-2}]\right)\\
&\quad+[q_1,\ldots, q_{s}]([q_{s+1},\ldots, q_{n-1}]q_n+[q_{s+1},\ldots, q_{n-2}]),\\
&=[q_1,\ldots, q_{s-1}]h[q_{s+2},\ldots, q_n]+
[q_1,\ldots, q_{s}][q_{s+1},\ldots, q_n].
\end{align*} 
	\end{proof}

\begin{proposition}[Property P--\ref{zigzagbis}]\label{prop:P4} Let $n,h,s$ be integers  such that $1\le s$, and let $q_1,\ldots, q_n$ be elements of a ring $R$. If there is a unit $u$ in $R$ commuting with each $q_i$, then
\[
[u ^{-1}q_1,uq_2,\ldots, u  ^{(-1)^n}q_n;h,s]=\begin{cases}
[q_1,\ldots,q_n;h,s], &\text{for even $n$;}\\
u^{-1}[q_1,\ldots,q_n;h,s], &\text{for odd $n$.}
\end{cases}
\]
\end{proposition}	

\begin{proof} Fix $s\ge 1$ and $h$ and proceed by induction on $n$. If $n=1,2$ then the property follows from Definition \ref{def:gCont}. 

If $s\ge n$ then $[u ^{-1}q_1,uq_2,\ldots, u  ^{(-1)^n}q_n;h,s]$ becomes $[u ^{-1}q_1,uq_2,\ldots, u  ^{(-1)^n}q_n]$ and the result follows from induction by using 
\begin{align*}
	[u ^{-1}q_1,uq_2,\ldots, u  ^{(-1)^n}q_n]&=[u ^{-1}q_1,uq_2,\ldots, u  ^{(-1)^{n-1}}q_{n-1}]u  ^{(-1)^n}q_n\\
	&\quad +[u ^{-1}q_1,uq_2,\ldots, u  ^{(-1)^{n-2}}q_{n-2}].
\end{align*}

If $s= n-1$ then  we have 
\begin{align*}
	[u ^{-1}q_1,uq_2,\ldots, u  ^{(-1)^n}q_n]&=[u ^{-1}q_1,uq_2,\ldots, u  ^{(-1)^{n-1}}q_{n-1}]u  ^{(-1)^n}q_n\\
	&\quad +[u ^{-1}q_1,uq_2,\ldots, u  ^{(-1)^{n-2}}q_{n-2}]h,
\end{align*}
and the result follows from induction.

If $s<n-1$ then Definition \ref{def:gCont} gives that 
\begin{align*}
[u ^{-1}q_1,uq_2,\ldots, u  ^{(-1)^n}q_n;h,s]&=[u ^{-1}q_1,uq_2,\ldots, u  ^{(-1)^{n-1}}q_{n-1};h,s]u^{(-1)^n}q_n\\
&\quad+[u ^{-1}q_1,uq_2,\ldots, u  ^{(-1)^{n-2}}q_{n-2};h,s].
\end{align*}
Then the induction hypothesis gives that 
\begin{align*}
[u ^{-1}q_1,uq_2,\ldots, u  ^{(-1)^{n-1}}q_{n-1};h,s]u^{(-1)^n}q_n=u ^{-1}[q_1,q_2,\ldots, q_{n-1};h,s]uq_n,
\end{align*}
if $n$ is even; and it gives that
\begin{align*}
[u ^{-1}q_1,uq_2,\ldots, u  ^{(-1)^{n-1}}q_{n-1};h,s]u^{(-1)^n}q_n=[q_1,q_2,\ldots, q_{n-1};h,s]u^{-1}q_n,
\end{align*}
if $n$ is odd. 

Furthermore, the induction hypothesis gives the following.

\begin{align*}	
[u ^{-1}q_1,uq_2,\ldots, u  ^{(-1)^{n-2}}q_{n-2};h,s]&=\begin{cases}
[q_1,q_2,\ldots, q_{n-2};h,s],&\text{if $n$ is even;}\\
u ^{-1}[q_1,q_2,\ldots, q_{n-2};h,s],&\text{if $n$ is odd.}
\end{cases}
\end{align*}	
As a consequence, the result follows.
\end{proof}
	
\begin{proposition}[Property P--\ref{determibis}]
\label{prop:P5}
 Let $n,h,s$ be integers  such that $1\le s$, and let $q_1,\ldots, q_n$ be elements of a commutative ring $R$. Then the generalized continuant $[q_1,\ldots, q_n;h,s]$
is the determinant of the tridiagonal $n\times n$ matrix $A=(a_{ij})$ with $a_{i,i}=q_i$
for $1\le i\le n$, $a_{i,i+1}=1$ for $1\le i<n$, $a_{s+1,s}=-h$ and $a_{i+1,i}=-1$ for $1\le i<n$
and $i\ne s$.
\end{proposition}
\begin{proof}
The result follows from using the Laplace expansion on the determinant along the last row.
\end{proof}

\begin{proposition}[Property P--\ref{reverbis}]
\label{prop:P6}
Let $n,h,s$ be integers  such that $1\le s$, and let $q_1,\ldots, q_n$ be elements of a commutative ring $R$. Then
 $[q_1,q_2,\ldots,q_n;h,n-s]=[q_n,\ldots,q_2,q_1;h, s]$. 
\end{proposition}
\begin{proof}
Apply
Property P--\ref{cuttingbis} on both sides of the equality, and then use Property P--\ref{PPreverse}.
\end{proof}

\providecommand{\bysame}{\leavevmode\hbox to3em{\hrulefill}\thinspace}
\providecommand{\MR}{\relax\ifhmode\unskip\space\fi MR }
\providecommand{\MRhref}[2]{
  \href{http://www.ams.org/mathscinet-getitem?mr=#1}{#2}
}
\providecommand{\href}[2]{#2}


\end{document}